
\input amstex
\documentstyle{amsppt}
\magnification=1200
\parskip 6pt
\NoBlackBoxes
\pagewidth{5.4in}
\hfuzz=5pt

\def\sR{\hbox{I\kern-.1667em\hbox{R}}}
\def\R{\hbox{I\kern-.1667em\hbox{R}}}
\def\N{\hbox{I\kern-.1667em\hbox{N}}}

\def\rn{{\R}^n}
\def\pair#1#2{\langle #1, #2 \rangle}

\def\vol{\hbox{vol}}
\def\spec{\hbox{spec}}
\def\hca{\hbox{hca}}
\def\vp{\hbox{vp}}
\def\prob{\hbox{P}}
\def\expec{\hbox{E}}
\def\mspec{\hbox{mspec}}

\def\b{{\Cal B}}


\def\Del#1#2{\frac{\partial #1}{\partial #2}}

\def\pair#1#2{\langle #1,#2 \rangle}

\def\chone{\hbox{1}}

\topmatter

\title
Dirichlet Spectrum and Heat Content
\endtitle

\author
Patrick McDonald and Robert Meyers
\endauthor

\affil
New College of Florida and The Courant
Institute of Mathematical Sciences
\endaffil

\address
New College of Florida and The Courant Institute
\endaddress

\email  ptm$\@$virtu.sar.usf.edu,  meyersr$\@$cims.nyu.edu
\endemail

\abstract
Let $M$ be a complete Riemannian manifold and $D\subset M$ a smoothly
bounded domain with compact closure.  We use Brownian motion and the
classic results on the Stieltjes moment problem to study the
relationship between the Dirichlet spectrum of $D$ and the heat
content asymptotics of $D.$  Central to our investigation is a
sequence of invariants associated to $D$ defined using exit time
moments.  We prove that our invariants determine that part of the
spectrum corresponding to eigenspaces which are not orthogonal to
constant functions, that our invariants  determine the heat content
asymptotics associated to the manifold, and that when the manifold is
a generic domain in Euclidean space, the invariants determine the
Dirichlet spectrum. 

\endabstract

\date March 24, 2002\enddate

\endtopmatter

\document

\heading
1:  Introduction
\endheading

Let $(M,g)$ be a complete Riemannian manifold and suppose
that $D\subset M$ is a smoothly bounded domain with compact
closure. Let $\spec(D)$ be the spectrum of the Laplace operator acting
on functions with Dirichlet boundary conditions.  We take the Laplacian
to be positive with elements of the spectrum listed in increasing
order with multiplicity.  We study the relationship between the
Dirichlet spectrum of $D$ and the heat content asymptotics of $D.$  We
recall the required facts:   

Let $p_D(x,y,t)$ be the heat kernel associated to $D,$ let $dg$ be the
volume form associated to the metric, and let 
$$\align
u(x,t) & = \int_D p_D(x,y,t) dg(y) \tag{1.1}
\endalign$$
be the solution to the initial value problem 
$$\aligned
\frac12 \Delta u & = \Del{u}{t} \hbox{ on } D \times (0,\infty) \\
u(x,0) & = \cases
                   1 & \hbox{ if } x \in D \\
                   0 & \hbox{ if } x \in \partial D  \endcases \\
u(x,t) & = 0 \hbox{ if } x \in \partial D
\endaligned\tag{1.2}$$
Let $q(t)$ be the heat content of $D$ at time $t:$
$$\align
q(t) & = \int_D u(x,t) dg. \tag{1.3}
\endalign$$
It is a theorem of van den Berg and Gilkey \cite{BG} that $q(t)$
admits a small time asymptotic expansion: 
$$\align
q(t) & \simeq \sum_{n=0}^\infty q_n t^{n/2} \tag{1.4}
\endalign$$
where the coefficients $q_n$ are locally computable geometric invariants
of $D$ (cf \cite{G} for a recent survey of results concerning heat
content).  We will refer to the coefficients occuring on right hand
side of \rom{(1.4)} as the {\it heat content asymptotics of $D$} and we
write 
$$\align
\hca(D) & = \{q_n\}_{n=0}^\infty. \tag{1.5}
\endalign$$
We note that, in contrast to the heat trace asymptotics, the heat
content asymptotics are {\it not} spectral.  Our results involve
relationships between the sets $\spec(D)$ and $\hca(D)$ for arbitrary
complete Riemannian manifolds $(M,g)$ and arbitrary smoothly bounded
domains with compact closure.  These  results arise naturally in the
context of probability and help to shed light on a wide range of
phenomena (estimates of the principal eigenvalue, comparison
theorems, isoperimetric phenomena, etc) tying probabilty to geometric
analysis.  The probabilistic tools involved are of two types: those
which give a probabilistic representation of the solution of boundary
value problems, and those involving the representation of nonnegative
sequences of real numbers as moments associated to a distribution
function (ie classical moment problems). To concisely state our
results, we recall the necessary material:  

Let $X_t$ be Brownian motion on $M.$  Let $\prob^x, \ x\in
M,$ be the family of probability measures charging
Brownian paths beginning at $x,$ and let $\expec^x$ be the
corresponding collection of expectation operators.  Let $\tau $ be the
first exit time of $X_t$ from $D:$ $$\tau = \inf\{t\geq 0: X_t \notin
D\}.$$ 
Then $u(x,t)$ defined in \rom{(1.2)} can be written as 
$$\align
u(x,t) & = \prob^x(\tau > t) \tag{1.6}
\endalign$$
Given \rom{(1.6)}, it is natural to consider the Laplace transform of
the random variable $\tau$ which is determined by the exit time
moments of $\tau.$  Thus, for $k$ a nonnegative integer, we are led to
consider the following nonnegative sequence of real numbers: 
$$\align
A_n & = \int_D \expec^x[\tau^n] dg. \tag{1.7}
\endalign$$
We write 
$$\align
\mspec(D) & = \{A_n\}_{n=0}^\infty \tag{1.8}
\endalign$$
and we note that $\mspec(D)$ is invariant under the action of the
isometry group of $M$ (cf 2.12-2.14).  Our first result is the following

\proclaim{Theorem 1.1} Let $(M,g)$ be a complete Riemannian
manifold, $D\subset M$ a smoothly bounded domain with compact closure.
For $\lambda \in \spec(D),$ let ${\Cal E}_\lambda(1)$ be orthogonal
projection of the constant function $1$ onto the eigenspace
correpsonding to $\lambda.$  Define constants $a_\lambda^2 \in \R $
by 
$$\align
a_\lambda^2 & = \int_D |{\Cal E}_\lambda(1)|^2 dg. \tag{1.9}
\endalign$$
Let $\spec^*(M)$ be the set whose elements are defined by 
$$\align
\spec^*(M) & = \left\{\lambda \in \spec(M) :  a_\lambda^2 \neq
0\right\}. \tag{1.10} 
\endalign$$
Then 
$$\align
\mspec(D) = \mspec(D') & {\text{ implies }} \spec^*(D) = \spec^*(D')
\endalign$$
and we say that $\mspec(D) $ determines $\spec^*(D).$   
\endproclaim

We remark that $\spec^*(D)$ is a set; in particular, it contains no
information concerning multiplicities.  

To prove Theorem 1.1 we note that the Stieltjes moment problem defined
by the sequence $\mspec(D)$ fixes a measure which determines both the
set $\spec^*(D)$ and the constants $a_\lambda^2$ defined by
\rom{(1.9)}.  This information, coupled to Theorem \rom{1.1},
determine the heat content asymptotics:  

\proclaim{Theorem 1.2}  Let $(M,g)$ be a complete Riemannian
manifold, $D\subset M$ a smoothly bounded domain with compact closure.
Then $\mspec(D)$ determines $\hca(D).$
\endproclaim

From the proof of Theorem 1.1 and Theorem 1.2, we obtain as a
corollary the fact that the information contained in $\spec^*(D)$ and
the {\it partition of volume} $\{a_\lambda^2\}_{\lambda \in \spec*(D)}$
determines $\hca(D)$ (cf Corollary 3.2).

We remark that results analogous to Theorem 1.1 and Theorem 1.2 hold
in the category of graphs and graph Laplacians (cf \cite{MM1}).  In
this context, there arises a natural Dirichlet series whose values at
positive integers gives the analog of invariants defined in
\rom{(1.1)} and whose values at negative integers gives the analog of
the heat content asymptotics of the associated graph domain (cf
\cite{MM1}).  In the context of domains in complete manifolds the same
is true; the relevant series arises as the Mellin transform of the
heat content which admits a meromorphic extension to the plane.  We
investigate properties of the meromorphic extension of this Dirichlet
series, characterizing the connection between special values,
residues, and the invariants of interest (cf Proposition
\rom{(2.1)}).  

Theorem 1.1, Theorem 1.2 and related results suggest that control of
the moment spectrum may be useful in studying a variety of geometric
phenomena including isoperimetric conditions (cf \cite{BS} and
\cite{M1} for related results, as well as the survey \cite{M2}), and
estimates for higher eigenvalues and spectral gaps.  In addition,
in the category of weighted graphs and discrete Laplacians, the moment
spectrum and heat content asymptotics distinguish analogues of the
isospectral nonisometric planar polygons of \cite{BCDS} (cf
\cite{MM2}), thus suggesting that heat content and Dirichlet spectrum
may provide a good collection of invariants for classifying smoothly
bounded domains up to isometry.

From the proof of Theorem 1.1 and Proposition 2.1 it is clear that
$\mspec(M)$ contains no information concerning multiplicity, nor does
it contain information concerning modes orthogonal to
constants.  Thus, in the presence of symmetry we expect that
$\mspec(M)$ will not provide full information concerning the Dirichlet
spectrum of the underlying domain.  Our final result indicates that
this occurrence is ``unusual''  when the manifold is a smoothly
bounded domain in Euclidean space with compact closure.

Recall, the collection of smoothly bounded domains in Euclidean space
with compact closure is naturally a Banach manifold.  We recall that a
property is {\it generic} for a Banach manifold if it holds for a set
of second category (ie, it holds for the complement of a countable
union of nowhere dense sets).  We prove:

\proclaim{Theorem 1.3}  For generic domains $D$ in $\rn,\ n \geq
2,$ $\mspec(D) $ determines $\spec(D).$
\endproclaim

Theorem 1.3 generalizes to domains in Riemannian manifolds of
dimension at least two.

\heading
2: Mellin transforms and Dirichlet series
\endheading

As in the introduction, let $(M,g)$ be a complete Riemannian
manifold and suppose that $D\subset M$ is a smoothly bounded domain
with compact closure.  Given continuous functions, $f, \ h,$ on $D,$
denote the natural pairing by $$\pair{f}{h} = \int_D fh dg.$$  Let
$\Delta$ be the Laplace operator and suppose that $\spec(D)$ is the
Dirichlet spectrum associated to $D.$  
\proclaim{Definition 2.1} Given $\lambda\in \spec(D),$ let ${\Cal
E}_\lambda(1)$ be orthogonal projection of the constant function $1$
onto the eigenspace associated to $\lambda.$  Let $a_\lambda^2$ be the
nonnegative real number defined by  
$$\align
a_\lambda^2 & =   \pair{{\Cal E}_\lambda(1)}{{\Cal
E}_\lambda(1)}. \tag{2.1}
\endalign$$
We call the set whose elements are given by $a_\lambda^2$ as
$\lambda $ runs through $\spec(D)$ a spectral partition of
volume and we write 
$$\align
\vp(D) & = \{a_\lambda^2\}_{\lambda \in \spec(D)}.
\endalign$$
\endproclaim
As suggested by Definition 2.1, the set $\vp(D)$ partitions the
volume amongst eigenspaces.  In particular 
$$\align
\sum_{\lambda\in \spec^*(D)} a_\lambda^2 & = \vol(D) \tag{2.2}
\endalign$$
where $\spec^*(D)$ is as in \rom{(1.10)}.  For $s$ complex,
$\hbox{Re}(s) \geq 0,$ we define 
$$\aligned
\zeta_D(s)  & = \sum_{\lambda \in \spec^*(D)} a_{\lambda}^2
 \left(\frac{2}{\lambda}\right)^s.  
\endaligned\tag{2.3}$$
We show that $\zeta_D(s)$ is closely related to the heat content
$q(t)$ defined in \rom{(1.3)}.  

Starting with the heat kernel written in terms of the spectral data,
we have 
$$\align
q(t) & = \sum_{\lambda \in \spec^*(D)} a_\lambda^2
e^{-\frac{\lambda t}{2}}.\tag{2.4} 
\endalign$$
We note that $q(t)$ is continuous and bounded on $[0,\infty).$  As
mentioned in the introduction, $q(t)$ admits a small time asymptotic
expansion given in \rom{(1.4)}.  For complex $s,$ $\hbox{Re}(s) >0,$
the Mellin transform of $q(t)$ is defined by
$$\align
{\Cal M}Q(s) & = \int_0^\infty q(t) t^s \frac{dt}{t}. \tag{2.5}
\endalign$$
Using \rom{(2.5)} we see that for $\hbox{Re}(s) >0,$
$$\align
{\Cal M}Q(s) & = \Gamma(s) \zeta_M(s) \tag{2.6}
\endalign$$
where $\zeta_D(s)$ is given by \rom{(2.3)} and $\Gamma(s)$ is the
gamma function.  By the standard theory of regularized series (cf
\cite{JL}), ${\Cal M}Q(s)$ admits a meromorphic extension to the plane
with poles restricted to lie at the negative half-integers.  In
addition, the poles are simple with residues given by
$$\align
\hbox{Residue}|_{s = -\frac{N}{2}} {\Cal M}Q(s) &  = q_N \tag{2.7}
\endalign$$
where $q_N$ is as in \rom{(1.4)} (ie the residues are given by the
heat content asymptotics).  This proves the first part of 

\proclaim{Proposition 2.1} For $\hbox{Re}(s)
>0,$ let $\zeta_D(s)$ be defined as in \rom{(2.3).}  Then
$\Gamma(s)\zeta_D(s)$ extends meromorphically to the complex plane
with poles restricted to lie at the negative half integers.  In
addition, the poles are simple and for $N$ a natural number, 
$$\align
\Gamma(N)\zeta_M(N) &  = \frac{1}{N} A_N \tag{2.8} \\
\hbox{Residue}|_{s =- \frac{N}{2}} \Gamma(s)\zeta_M(s) &  = q_N \tag{2.9}
\endalign$$
where $A_N$ is given by \rom{(1.7)} and $q_N$ is given by
\rom{(1.4)}. 
\endproclaim

\demo{Proof}  The claim \rom{(2.9)} follows immediately from \rom{(2.6)}
and \rom{(2.7)}.  To see that \rom{(2.8)} holds, let $X_t$ be Brownian
motion on $M,$ and let $\tau$ be the first exit time from $D.$
Let $$h(x,s) = \expec^x[e^{-s\tau}]. \tag{2.10}$$ Then $h$ is the unique
solution of the Dirichlet problem  
$$\aligned
\frac12 \Delta h - sh & = 0 \hbox{ on } D \times (0,\infty) \\
h & = 1 \hbox{ on }  D \times \{0\}.
\endaligned\tag{2.11}$$
Expanding $h(x,s)$ using power series and using \rom{(2.11),} we see
that the exit time moments can be defined by recursive solution of
Poisson problems.  More precisely, suppose that  
$$\aligned
\frac12 \Delta u_1 + 1 & = 0 \hbox{ on } D \\
u_1 & = 0 \hbox{ on }\partial D
\endaligned\tag{2.12}$$
and
$$\aligned
\frac12 \Delta u_k + ku_{k-1} & = 0 \hbox{ on } D \\
u_k & = 0 \hbox{ on } \partial D.
\endaligned\tag{2.13}$$
Then $$\expec^x[\tau^k] = u_k(x).\tag{2.14}$$  Thus, if $\phi_\lambda$ is a
normalized eigenfunction, we have
$$\align
\pair{u_k}{\phi_\lambda} &= -\frac{2}{\lambda}
\pair{u_k}{\frac12 \Delta \phi_\lambda} \\
 & = \frac{2}{\lambda}k \pair{u_{k-1}}{\phi_\lambda}.\tag{2.15}
\endalign$$
Writing $$1 = \sum_{\lambda \in \spec(M)}
\pair{1}{\phi_\lambda}\phi_\lambda,$$ the proposition follows.
\enddemo

\heading
3: Proof of Theorem 1.1 and Theorem 1.2
\endheading

We begin with a corollary of Proposition \rom{2.1}:

\proclaim{Corollary 3.1} For $A_n$ as defined in \rom{(1.7)}, set
$$\mu_n = \frac{A_n}{n!}. \tag{3.1}$$  Then the collection $\{\mu_n
\}$ satisfies Carleman's condition:
$$\align
\sum \mu_{2n}^{-\frac{1}{2n}} & = \infty. \tag{3.2}
\endalign$$
\endproclaim

\demo{Proof} From Proposition \rom{(2.1)}, $\mu_n = \zeta_D(n).$
Thus,
$$\align
\mu_n & = \sum_{\lambda \in \spec^*(D)} a_\lambda^2
\left(\frac{2}{\lambda} \right)^n \\ 
 & \leq \left(\frac{2}{\lambda_1}\right)^n \hbox{vol}(D).
\endalign$$
\enddemo
\demo{Proof of Theorem 1.1} Let $\mu_n$ be as defined in
\rom{(3.1)} and note that $\mu_n >0$ for all nonnegative integers $n.$
Define a bounded, nondecreasing function, $\psi : [0,\infty ) \to
[0,\infty)$ by 
$$\align
\psi(x) & = \sum_{\lambda \in \spec^*(M)} a_\lambda^2
\chone_{[\frac{1}{\lambda}, \infty)}(x)\tag{3.3}
\endalign$$
where $\chone_{[\frac{1}{\lambda}, \infty)}(x)$ is the indicator
function of the interval $[\frac{1}{\lambda},\infty).$  Then $\psi$
solves the Stieltjes moment problem for the moments $\mu_n:$
$$\align
\mu_n & = \int_0^\infty x^n d\psi .
\endalign$$
Recall the classic result of Carleman (cf \cite{A}): 
\proclaim{Theorem} Suppose $\{\mu_n\}$ is a sequence of nonnegative
real numbers.  If \rom{(3.2)} holds, then the Stieltjes moment problem
for the sequence $\{\mu_n\}$ is determined.
\endproclaim
By Corollary \rom{3.1}, the moments satisfy Carleman's condition
and thus the unique solution of the Stieltjes Moment Problem is given
by \rom{(3.3)}.  Thus, the sequence $\{A_n\}$ determines both the set
$\spec^*(M)$ (the discontinuities of  $\psi(x)$), as well as the
collection of jumps, $\vp(D).$  This proves Theorem
1.1. 

\enddemo

From the proof of Theorem 1.1 we immediately conclude:

\proclaim{Corollary 3.1} Let $D\subset M $ be a smoothly bounded
domain with  compact closure.  Then $\mspec(D)$ determines $\vp(D).$  
\endproclaim

\demo{Proof of Theorem 1.2} Theorem 1.2 follows immediately from
Theorem 1.1, Corollary 3.1, and \rom{(2.4)}.  
\enddemo

Finally, we give a relationship between Dirichlet spectrum and heat
content asymptotics.  

\proclaim{Corollary 3.2} Let $D\subset M $ be a smoothly bounded
domain with  compact closure.  Then $\spec^*(D)\cup \vp(D)$ determines
$\hca(D).$   
\endproclaim

\demo{Proof}  This is immediate from the definition of $\zeta_D(s)$
and Proposition 2.1.  In fact, from Corollary 3.1 and \rom{(2.4)} it
is clear that the heat content (not just the asymptotics) is
determined. 
\enddemo

\heading
4: Proof of Theorem 1.3
\endheading
In this section we consider $C^k$-domains with compact closure in
$\rn, \ n \geq 2.$

Let $k> n+2$ and let $\b$ be the collection of $C^k$-domains in
$\rn$ with compact closure.  Recall, $\b$ is a Banach manifold:
Given $b \in \b,$ we identify $b$ with its boundary, $\partial b.$
The tubular neighborhood theorem identifies a neighborhood of
$\partial b$ in $\rn$ with sections of the normal bundle to the
boundary of $b,$ denoted $C^k(\partial b, N\partial b).$  Pairing
with the outward pointing unit normal vector gives an isomorphism
between $C^k(\partial b, N\partial b)$ and $C^k(\partial b).$  We
identify domains near $b$ by identifying their boundaries as those
obtained by flow in the normal direction prescribed by elements of
$C^k(\partial b).$  More precisely, if $\nu$ is the outward
pointing unit normal vector along the boundary of $b$ and $f\in
C^k(\partial b),$ then for $\epsilon$ small enough, the set
$$\partial b_\epsilon = \{y \in \rn: y =\sigma + \epsilon
f(\sigma) \nu(\sigma), \ \ \sigma \in \partial
b\}\tag{4.1}$$bounds a $C^k$-domain in $\rn$ and $\epsilon \to
b_\epsilon$ where $b_\epsilon$ is the domain bounded by $\partial
b_\epsilon$ is a smooth curve in $\b$ passing through $b$ at
$\epsilon = 0.$  This provides an identification of a neighborhood
of $b \in \b$ with a neighborhood of $0$ in $C^k(\partial b),$
which shows that $\b$ is a Banach manifold.  In addition, the
construction indicates that that there is a natural choice for the
tangent space of $b\in \b:$
$$\align
T_b\b & \simeq C^k(\partial b). \tag{4.2}
\endalign$$

It is a theorem of Uhlenbeck (\cite{U} also \cite{CV}) that the collection of
$C^k$-domains $b \in \b$ for which all Dirichlet eigenvalues have
multiplicity one is open and dense in $\b.$  In the sequel, we adopt
Uhlenbeck's approach to establish a generic property useful for our purposes.
We begin by recalling the necessary notation.

Let $H_k(\rn)$ be the Sobolev space of functions on $\rn$ with
distributional derivatives up through order $k$ which are $L^2.$  Let
$H_{k,0}(\rn)\subset H_k(\rn)$ be the closure of the space of smooth
functions on $\rn,$ and let $S_k(\rn)$ be the unit ball in
$H_k(\rn).$  Let $\phi:S_k(\rn) \times \R \times {\Cal B} \to
H_{k-2}(\rn)$ be defined by $\phi(u,\lambda,b) = \Delta_bu - \lambda u$
where $\Delta_b$ is the Dirichlet Laplacian on $b.$  Let $Q\subset
S_k(\rn) \times \R \times \b$ be defined by $Q = \phi^{-1}(0).$  Then
$Q$ is the collection of domains, their Dirichlet spectrum, and their 
corresponding normalized Dirichlet eigenfunctions.  It is a corollary of the
Sard-Smale theorem that $Q$ is a Banach submanifold of $ H_k(\rn)
\times \R \times \b.$  The tangent space at a point $(u,\lambda,b) \in
Q$ is given by
$$\aligned
T_{(u,\lambda,b)} & \simeq \left\{ (v,\eta, f) \in H_{k,0}(b) \times \R
\times T_b\b: \int_b uv = 0 , \right. \\
 &  \ \hbox{   } \ \left.  (\Delta_b +\lambda)v + \eta u +D_b\phi(f) =
 0\right\}  
\endaligned\tag{4.3}$$
where $D_b\phi$ is the derivative of the function $\phi$ with respect
to ${\Cal B}.$ 

\proclaim{Definition} Let $b\in \b.$  We say that $b$ has Property M if
for all $\lambda \in \spec(b),$
$$\align
\pair{{\Cal E}_\lambda(1)}{{\Cal E}_\lambda(1)}& \neq 0 \tag{4.4}
\endalign$$
where ${\Cal E}_\lambda$ is projection on the eigenspace corresponding
to $\lambda.$
\endproclaim

To see that Property M is generic, we define a function $I:Q \to \R$
by
$$\align
I(u,\lambda,b) & = \int_b u(x) dx \tag{4.5}
\endalign$$
where $dx$ denotes Lebesgue measure.  We note that $I$ is clearly
$C^k,$ and thus $I^{-1}(\rn \setminus \{0\})$ is open.  We will show
that $D_bI,$ the derivative of $I$ with respect to domain variations, is
always surjective. 

To see that this is the case fix $(u,\lambda, b)\in Q.$  An
infinitesimal variation of the domain $b$ is given by fixing an
element $f \in C^k(\partial b).$  We denote by $\delta u$ the
corresponding infinitesimal change in $u$ and by $\delta \lambda$
the corresponding infinitesimal change in $\lambda.$  A
straightforward computation then gives:
$$\align
D_bI(\delta u ,\delta \lambda, f) & = \int_b \delta u(x) dx. \tag{4.6}
\endalign$$
Using Hadamard's classic results on the variation of Green's functions
for perturbed domains (cf \cite{H}, \cite{GS}), we have an expression
for $\delta u(x):$ 
$$\align
\delta u(x) & = -\int_{\partial b} f(\sigma) \Del{u}{\nu}(\sigma)
\Del{G}{\nu}(x,\sigma) d\sigma \tag{4.7}
\endalign$$
where $d\sigma$ is the induced surface measure on the boundary and $G$
is the Green's function for $D.$  As pointed out by Uhlenbeck for a
similar computation, it is a corollary of unique continuation that $
\Del{u}{\nu}(\sigma) \Del{G}{\nu}(x,\sigma) $ is not identically zero.
We conclude that we can find $f$ such that $D_bI(\delta u ,\delta
\lambda, f) \neq 0.$  As a corollary, we obtain 

\proclaim{Theorem 4.1} Let $\b$ be the Banach manifold of $C^k$-domains
with compact closure.  Then Property $M$ is generic for $\b.$
\endproclaim

\demo{Proof of Theorem 1.3} From Theorem \rom{1.1}, we know that
$\spec^*(D)$ is determined by $\mspec(D).$  By Uhlenbeck's theorem,
for a dense open set of domains, all eigenvalues have multiplicity
one.  By Theorem 4.1, for domains with all eigenvalues of multiplicity
one, it is generically the case that $\spec^*(D) = \spec(D).$  This
concludes the proof of Theorem 1.3.
\enddemo

\enddocument